# ASYMPTOTICS OF SOLUTIONS TO SEMILINEAR STOCHASTIC WAVE EQUATIONS

By Pao-Liu Chow[1]

*Wayne State University*

Large-time asymptotic properties of solutions to a class of semilinear stochastic wave equations with damping in a bounded domain are considered. First an energy inequality and the exponential bound for a linear stochastic equation are established. Under appropriate conditions, the existence theorem for a unique global solution is given. Next the questions of bounded solutions and the exponential stability of an equilibrium solution, in mean-square and the almost sure sense, are studied. Then, under some sufficient conditions, the existence of a unique invariant measure is proved. Two examples are presented to illustrate some applications of the theorems.

**1. Introduction.** Semilinear stochastic wave equations arise as mathematical models to describe nonlinear vibration or wave propagation in a randomly excited continuous medium. To be specific, the equation may take the form

$$(1.1) \qquad \partial_t^2 u(x,t) = c^2 \Delta u - 2\alpha\, \partial_t u(x,t) + f(u) + \sigma(u)\dot{W}(x,t)$$

in a bounded domain $\mathcal{D}$ in $\mathbf{R}^d$, subject to some homogeneous boundary and initial conditions to be specified later. Here $\partial_t = \frac{\partial}{\partial t}$, $\Delta$ is the Laplacian operator, and $c$ and $2\alpha$ are some positive constants known as the wave speed and the damping coefficient, respectively. The nonlinear functions $f$ and $\sigma$ are given, and $\dot{W}(x,t) = \partial_t W(x,t)$ is a spatially dependent white noise, where $W(x,t)$ is a Wiener random field. In previous papers [3, 4], we studied the local and global solutions of this type of equation without damping ($\alpha = 0$), where the nonlinear terms $f$ and $\sigma$ may admit a polynomial growth. As

Received October 2004; revised July 2005.
[1]Supported in part by a grant from the Academy of Applied Science.

a sequel to our previous work, this paper is concerned with some qualitative asymptotic behavior of solutions to equation (1.1) in a bounded domain $\mathcal{D}$ as $t \to \infty$. In addition to the global existence of solutions, we are interested in the questions of boundedness, asymptotic stability and the existence of a stationary solution or an invariant measure. For a solution of the wave equation to reach a statistical equilibrium, it is imperative to include the damping term in equation (1.1) so that, in the physical term, the fluctuation–dissipation principle may hold. As a simple example, consider the randomly perturbed wave equation in one dimension,

$$(1.2) \quad \begin{aligned} \partial_t^2 u &= c^2 \partial_x^2 u - 2\alpha \, \partial_t u + \dot{W}(x,t), \qquad t > 0, \ x \in \mathcal{D} = (0,\pi), \\ u(x,0) &= h(x), \qquad \partial_t u(x,0) = 0; \qquad u(0,t) = u(\pi,t) = 0, \end{aligned}$$

where $h$ is a given continuous function and the Wiener field $W$ is assumed to have the Fourier series representation

$$W(x,t) = \sum_{n=1}^{\infty} \sigma_n b_n(t) \phi_n(x),$$

where $\{b_n(t)\}$ is a sequence of independent copies of standard Brownian motions in one dimension, $\{\sigma_n\}$ is a sequence of reals such that $\sum_{n=1}^{\infty} \sigma_n^2 < \infty$ and $\phi_n = \sqrt{2/\pi} \sin nx$, $n = 1, 2, \ldots$, are the normalized eigenfunctions associated with the problem (1.2). Then, by means of the eigenfunction expansion, (1.2) can be formally solved in the case $c > \alpha$ to give

$$(1.3) \qquad u(x,t) = \sum_{n=1}^{\infty} u_n(t) \phi_n(x),$$

where

$$(1.4) \qquad u_n(t) = h_n e^{-\alpha t} \cos \omega_n t + \frac{\sigma_n}{\omega_n} \int_0^t e^{-\alpha(t-s)} \sin \omega_n(t-s) \, db_n(s)$$

with $h_n = \int_0^\pi h(x) \phi_n(x) \, dx$ and $\omega_n = \sqrt{(nc)^2 - \alpha^2}$ for $n = 1, 2, \ldots$. By some simple calculations, we obtain the mean

$$E u_n(t) = h_n e^{-\alpha t} \cos \omega_n t \to 0$$

and the variance

$$\mathrm{Var}\{u_n(t)\} = \left(\frac{\sigma_n}{\omega_n}\right)^2 \int_0^t e^{-2\alpha s} \sin^2(\omega_n s) \, ds \to \frac{1}{4\alpha} \left(\frac{\sigma_n}{nc}\right)^2$$

as $t \to \infty$. Then it follows from (1.4) that the solution $u(x,t)$ is a Gaussian random field with mean $Eu(x,t) \to 0$, and covariance function

$$\lim_{t \to \infty} \mathrm{Cov}\{u(x,t), u(y,t)\} = \sum_{n=1}^{\infty} \frac{1}{4\alpha} \left(\frac{\sigma_n}{nc}\right)^2 \phi_n(x) \phi_n(y).$$

In fact it can be shown that the solution $u(\cdot, t)$ converges in the mean-square to a Gaussian random field $\hat{u}(\cdot)$ with the above covariance function and its probability law is the invariant measure for equation (1.2). On the other hand, without the damping ($\alpha = 0$), we would have $Eu_n(t) = h_n \cos nct$ and

$$\text{Var}\{u_n(t)\} = (1/2)(\sigma_n/nc)^2[t - (1/2nc)\sin 2nct] \to \infty \quad \text{as } t \to \infty.$$

So the asymptotic solution will cease to exist. Clearly this can also happen in the nonlinear case. In fact it was shown that, with a cubic nonlinearity, the solution may explode in finite time [3] unless there exists a certain energy bound. As to be seen, the dissipation and the energy bound for a semilinear wave equation such as (1.1) are two major ingredients to ensure proper asymptotic behavior of its solutions.

Our initial work on semilinear stochastic wave equations [3] was stimulated by two interesting papers by Mueller [18, 19] on the existence of large-time solutions to some nonlinear heat and wave equations with noise. If such a solution exists, it is natural to investigate its asymptotic behavior as $t \to \infty$. By a semigroup approach, asymptotic solutions to semilinear stochastic evolution equations have been studied by many authors. For the problems of boundedness and stability, see, for example, the papers [2, 10, 12, 13], and for the existence of invariant measures, we mention the articles [5, 8, 14, 16], and the book [7] for further references. In concrete terms, most of the above-mentioned results are applicable to the parabolic or dissipative type of stochastic partial differential equations. The asymptotic solution of a stochastic hyperbolic or wavelike equation was studied in [15] by the method of averaging. To our knowledge the asymptotic solutions of the semilinear wave equations under consideration have not been treated in the literature. For the deterministic case, the analysis of hyperbolic equations relies heavily on the so-called energy method ([23], [25], Chapter 4). Therefore the associated energy function plays an important role in the asymptotic analysis. Similarly we shall adopt the stochastic version of the energy method in the current study. In fact, to obtain the crucial exponential estimates, it is necessary to introduce a pseudo energy function, which can be interpreted physically as adding an artificial damping to the system. For some related works on stochastic wave equations, we mention the interesting papers [17, 20, 21], among many others.

**2. Summary of results.** In Section 3 we present three technical lemmas. In Lemma 3.1 we prove the existence of a unique solution and the energy equation for a linear stochastic wave equation. By introducing a pseudo energy function, a key exponential estimate is established in Lemma 3.2. Then it is shown in Lemma 3.3 that the pseudo energy function is equivalent to the usual energy function.

The global solution to a class of semilinear stochastic wave equations of the form (4.1) is treated in Section 4. Under the locally bounded, local Lipschitz conditions in the Sobolev space $H^1$ and an energy inequality given by conditions (A1)–(A4), the results of the existence and uniqueness of a global solution are stated and proved in Theorem 4.1. The proof is based on a smooth $H^1$-truncation technique and some probabilistic inequalities.

In Section 5 we consider the boundedness of solutions in a mean-square sense as $t \to \infty$. Assuming that, in addition to Conditions A, the nonlinear terms satisfy a set of growth conditions (B1)–(B3), Theorem 5.1 shows that the solution is bounded in mean-square, while with a slightly stronger assumption, it is proved in Theorem 5.2 that the solution is ultimately bounded in the mean-square sense.

Then some questions of the asymptotic stability of the null solution are considered in Section 6. Under Conditions B with an exponential integrability condition on the parametric functions $\theta(t)$ and $\rho(t)$, Theorem 6.1 shows that the null solution is asymptotically, exponentially stable in mean-square. If $\theta = \rho \equiv 0$, as stated and proved in Theorem 6.3, the null solution becomes exponentially stable almost surely.

So far the stochastic wave equations under consideration admit that nonlinear terms satisfy only a local Lipschitz condition in the space $H^1$. In particular, the nonlinear terms are allowed to have a polynomial growth. For the existence of an invariant measure, this poses a challenging open problem as yet to be resolved. Even in the case of globally Lipschitzian nonlinearity, the existence result does not seem to have been proven. Hence, in Section 7, we shall prove an existence theorem (Theorem 7.1) by assuming that the nonlinear terms are globally Lipschitzian and have linear growth. Technically our proof follows the approach of Da Prato and Zabczyk [7] by adapting their method for the strongly dissipative equation, such as a parabolic equation, to our hyperbolic problem. Finally two examples are provided in Section 8 to illustrate some applications of our theorems.

**3. Energy equation and exponential estimate.** Let $\mathcal{D} \subset \mathbf{R}^d$ be a bounded domain with a smooth, say, $C^2$ boundary $\partial \mathcal{D}$. We set $H := L^2(\mathcal{D})$ with the inner product and norm denoted by $(\cdot, \cdot)$ and $\|\cdot\|$, respectively. Let $H^k = W^{k,2}(\mathcal{D})$ be the $L^2$ Sobolev space of order $k$ with norm $\|\cdot\|_k$, and denote by $H_0^1$ the closure in $H^1$ of the set of all $C^1$ functions with compact support in $\mathcal{D}$. The dual space of $H^1$ is given by $H^{-1}$ [1].

Let $(\Omega, \mathcal{F}, P)$ be a complete probability space for which a filtration $\mathcal{F}_t$ of sub-$\sigma$-fields of $\mathcal{F}$ is given. Let $W(x, t)$, $x \in \mathcal{D}$, $t \geq 0$, be a continuous Wiener random field defined in this space with $W(x, 0) = 0$. It has a zero mean, $EW(x, t) = 0$ and covariance

(3.1) $\quad E[W(x, t) W(y, s)] = (t \wedge s) r(x, y), \qquad x, y \in \mathcal{D},$

where $(t \wedge s) = \min(t,s)$ for $0 \leq t, s \leq T$ and the covariance function $r(x,y)$ is bounded so that

$$\sup_{x \in \mathcal{D}} r(x,x) \leq r_0. \tag{3.2}$$

Let $\sigma(x,t) = \sigma(x,t,\omega)$ for $t \geq 0, x \in \mathcal{D}$ and $\omega \in \Omega$ be a continuous $\mathcal{F}_t$-predictable random field satisfying the condition

$$E \int_0^T \|\sigma(\cdot,t)\|^p \, dt < \infty \tag{3.3}$$

for $p \geq 2$. Then it can be shown that the stochastic integral

$$M(x,t) = \int_0^t \sigma(x,s) W(x,ds), \qquad t > 0, \ x \in \mathcal{D}, \tag{3.4}$$

is well defined and $M_t = M(\cdot,t)$ is a continuous $H$-valued $\mathcal{F}_t$ martingale (see the Appendix). It has mean $EM(x,t) = 0$ and covariation operator $Q_t$ defined by

$$\langle (M.,g), (M.,h) \rangle_t = \int_0^t (Q_s g, h) \, ds$$

for any $g, h \in H$, where the kernel function $q(x,y,t)$ of $Q_t$, defined by

$$(Q_t g)(x) = \int_{\mathcal{D}} q(x,y,t) g(y) \, dy,$$

is given by

$$q(x,y,t) = r(x,y) \sigma(x,t) \sigma(y,t).$$

In view of conditions (3.2) and (3.3), it can be shown that (see the proof of Theorem A.1)

$$E\|M_t\|^p \leq C_p(T) E \int_0^T \|\sigma(\cdot,t)\|^p \, dt$$

for some positive constant $C_p(T)$. It is worth noting that the stochastic integration in (3.4) is taken with respect to an $L^p$-bounded integrand $\sigma_t$, instead of a Hilbert–Schmidt operator-valued process as usually done (see [6], Chapter 4). This version of stochastic integral will be needed later on to deal with equations with pointwise (in $x$) multiplicative noises (see Example 1). Since we have not been able to find a reference for this type of integral, it will be defined in the Appendix.

Now we consider the initial boundary value problem for the linear damped hyperbolic equation with a random perturbation,

$$[\partial_t^2 + 2\alpha \partial_t - A(x,D)] u(x,t) = f(x,t) + \partial_t M(x,t), \qquad 0 < t < T,$$
$$u(x,0) = u_0(x), \qquad \partial_t u(x,0) = v_0(x), \qquad x \in \mathcal{D}, \tag{3.5}$$
$$u(\cdot,t)|_{\partial \mathcal{D}} = 0,$$

where $\alpha$ is a positive parameter, $D = \partial_x$ denotes the gradient operator and $A(x, D)$ is a strongly elliptic operator of second order of the form

$$(3.6) \qquad A(x,D)\varphi(x) = \sum_{i,j=1}^{d} \partial_{x_i}[a^{ij}(x)\,\partial_{x_j}\varphi(x)] - b(x)\varphi(x).$$

In addition, the coefficients $a^{ij} = a^{ji}$ and $b$ are assumed to be smooth functions that satisfy

$$a_0(1+|\xi|^2) \le \sum_{i,j=1}^{d} a^{ij}(x)\xi_i\xi_j + b(x)|\xi|^2 \le a_1(1+|\xi|^2), \qquad \xi, x \in \mathcal{D},$$

for some constants $a_1 \ge a_0 > 0$.

To consider (3.5) as an Itô equation in a Hilbert space, we set $u_t = u(\cdot, t), v_t = v(\cdot, t)$ and so on, and rewrite it as

$$(3.7) \qquad \begin{aligned} du_t &= v_t\, dt, \\ dv_t &= [Au_t - 2\alpha v_t + f_t]\, dt + dM_t, \qquad 0 < t < T, \\ u_0 &= g, \qquad v_0 = h, \end{aligned}$$

where the domain $\mathcal{D}(A) = H^2 \cap H_0^1$, $g \in H^1, h \in H$ and $M_t$ is regarded as an $H$-valued Wiener martingale. Condition (3.6) implies that $(-A)$ is a self-adjoint, strictly positive linear operator in $H = L^2(D)$ and its square root $B = \sqrt{-A}$ is also a self-adjoint, strictly positive operator with domain $\mathcal{D}(B)$, which is a Hilbert space under the inner product $(g, h)_B := (Bg, Bh)$ and norm $\|g\|_B = (Bg, Bg)^{1/2}$ (see [24], Chapter 1). Since $D(B) \cong H_1$, for convenience, we define $\|\cdot\|_1 = \|\cdot\|_B$ in the subsequent analysis. As usual, the Itô differential equation (3.7) is interpreted as a stochastic integral equation:

$$(3.8) \qquad \begin{aligned} u_t &= u_0 + \int_0^t v_s\, ds, \\ v_t &= v_0 + \int_0^t Au_s\, ds - 2\alpha \int_0^t v_s\, ds + \int_0^t f_s\, ds + M_t. \end{aligned}$$

Introduce the Hilbert space $\mathcal{H} = (H^1 \times H)$ with $\mathcal{H}_0 = (H_0^1 \times H)$, equipped with the norm defined by

$$\|\phi\|_{\mathcal{H}} = \{\|u\|_1^2 + \|v\|^2\}^{1/2} = \{\|Bu\|^2 + \|v\|^2\}^{1/2}$$

for any $\phi = (u; v) \in \mathcal{H}$. Let $\mathcal{H}^* = (H^{-1} \times H)$ denote the dual space of $\mathcal{H}$. Define the *energy function* $\mathbf{e}(\cdot) \colon \mathcal{H} \to \mathbf{R}^+ = [0, \infty)$ as

$$(3.9) \qquad \mathbf{e}(\phi) := \mathbf{e}(u; v) = \|Bu\|^2 + \|v\|^2 \qquad \text{for } \phi = (u; v) \in H^1 \times H.$$

Notice that the norm $\|\phi\|_{\mathcal{H}} = \sqrt{\mathbf{e}(\phi)}$ is also called an energy norm. In what follows, we denote the $\mathcal{H}$-norm $\|\cdot\|_{\mathcal{H}}$ simply by $\|\cdot\|$ when there is no confusion.

Now by regarding (3.8) as a stochastic evolution equation in $\mathcal{H}^*$ in the distributional sense, we have the following lemma:

LEMMA 3.1 (Energy equation). *For $\phi_0 = (u_0; v_0) \in \mathcal{H}$, let $f_t$ be a continuous predictable process in $H$ and let $M_t$ be a continuous $H$-valued martingale with covariation operator $Q_t$ such that*

$$E\left\{\int_0^T \|f_t\|^2\, dt + \int_0^T \mathrm{Tr}\, Q_t\, dt\right\} < \infty, \tag{3.10}$$

*where* $\mathrm{Tr}$ *denotes the trace operator in* $\mathbf{H}$. *Then equation (3.8) or (3.7) has a unique solution $\phi_t = (u_t; v_t)$ which is a continuous $\mathcal{H}$-valued semimartingale. Moreover, it satisfies the energy equation*

$$\begin{aligned}\mathbf{e}(\phi_t) &= \mathbf{e}(\phi_0) - 4\alpha \int_0^t \|v_s\|^2\, ds + 2\int_0^t (v_s, f_s)\, ds \\ &\quad + 2\int_0^t (v_s, dM_s) + \int_0^t \mathrm{Tr}\, Q_s\, ds \qquad a.s.\end{aligned} \tag{3.11}$$

*for $t \in [0, T]$, where the energy function $\mathbf{e}(\cdot)$ on $\mathcal{H}$ is defined by (3.9). Moreover, the inequality*

$$E \sup_{t \leq T} \mathbf{e}(\phi_t) \leq C_1 + C_2 E \int_0^T \{\|f_s\|^2 + \mathrm{Tr}\, Q_s\}\, ds \tag{3.12}$$

*holds, where the constants $C_1, C_2$ depend on $p, T$ and the initial conditions.*

PROOF. Since the idea of the proof is similar to that of Lemma 2.1 in [4] with the Laplacian replaced by $A$, we will only sketch the proof. The only difference is that, instead of using Friedrichs' mollifying approximation, we adopt a finite-dimensional projection.

To this end, since $A$ is strongly elliptic and self-adjoint, it has a complete orthonormal set of eigenfunctions $\{\varphi_n\}$ with corresponding eigenvalues $\{\lambda_n\}$. Let $P_n: H \to H_n$ be defined by $P_n h = \sum_{k=1}^n (h, \varphi_k)\varphi_k$, where $H_n$ is a finite-dimensional subspace of $H^2$ spanned by $\{\varphi_1, \ldots, \varphi_n\}$. Apply the projector $P_n$ to equation (3.7) to get

$$\begin{aligned}du_t^n &= v_t^n\, dt, \\ dv_t^n &= [Au_t^n - 2\alpha v_t^n + f_t^n]\, dt + dM_t^n, \qquad 0 < t < T, \\ u_0^n &= g^n, \qquad v_0^n = h^n,\end{aligned} \tag{3.13}$$

where we set $u_t^n = P_n u_t, \ldots, h^n = P_n h$. The finite-dimensional linear system has a unique $\mathcal{F}_t$-adapted continuous solution $(u_t^n; v_t^n)$ in $(H_n \times H_n) \subset (H_0^1 \times$

$H$). In particular, by the Itô formula, the following energy equation holds:

$$\mathbf{e}(u_t^n, v_t^n) = \mathbf{e}(g^n, h^n) - 4\alpha \int_0^t \|v_s^n\|^2 \, ds + 2 \int_0^t (v_s^n, f_s^n) \, ds$$
(3.14)
$$+ 2 \int_0^t (v_s^n, dM_s^n) + \int_0^t \operatorname{Tr} Q_s^n \, ds \quad \text{a.s.}$$

By means of the simple inequality $2(v^n, f^n) \leq (\alpha \|v^n\|^2 + \frac{1}{\alpha} \|f^n\|^2)$ and the B–D–G (Burkholder–Davis–Gundy) inequality ([6], page 82),

$$E \sup_{0 \leq t \leq T} \left| 2 \int_0^t (v_s^n, dM_s^n) \right| \leq 8 E \left\{ \int_0^T \|v_s^n\|^2 \operatorname{Tr} Q_s^n \, ds \right\}^{1/2}$$

$$\leq 8 E \left\{ \sup_{0 \leq t \leq T} \|v_t^n\| \left| \int_0^T \|v_s^n\| \operatorname{Tr} Q_s^n \, ds \right|^{1/2} \right\}$$

$$\leq \tfrac{1}{2} E \sup_{0 \leq t \leq T} \|v_t^n\|^2 + 36 \int_0^T \operatorname{Tr} Q_s^n \, ds,$$

we can deduce from (3.14) that

(3.15) $$E \sup_{0 \leq t \leq T} \mathbf{e}(u_t^n, v_t^n) \leq C_1 E \int_0^T \{\|f_s\|^2 + \operatorname{Tr} Q_s^n\} \, ds$$

for some constant $C_1 > 0$. Let $X = L^2(\Omega; C([0;T], H_0^1 \times H))$ with norm $\|(u;v)\|_X = \{E \sup_{0 \leq t \leq T} \mathbf{e}(u_t, v_t)\}^{1/2}$. Then $X$ is known to be a separable, reflexive Banach space ([22], page 218). In view of (3.15), the sequence $\{(u^n; v^n)\}$ is bounded in $X$ so that there exists a subsequence $\{(u^{n_k}; v^{n_k})\}$ that converges weakly to $(u; v) \in X$.

In fact, we can show that the subsequence converges strongly in $X$. To do so, denote the subsequence again by $\{(u^n; v^n)\}$ and set $(u^{mn}; v^{mn}) = (u^m; v^m) - (u^n; v^n)$. It suffices to show that $\{(u^n; v^n)\}$ is a Cauchy sequence in $X$ so that $\|(u^{mn}; v^{mn})\|_X \to 0$ as $m, n \to \infty$. In view of (3.13) and (3.14), the difference sequence satisfies the energy equation

$$\mathbf{e}(u_t^{mn}; v_t^{mn}) = \mathbf{e}(g^{mn}; h^{mn}) - 4\alpha \int_0^t \|v_s^{mn}\|^2 \, ds + 2 \int_0^t (v_s^{mn}, f_s^{mn}) \, ds$$
(3.16)
$$+ 2 \int_0^t (v_s^{mn}, dM_s^{mn}) + \int_0^t \operatorname{Tr} Q_s^{mn} \, ds \quad \text{a.s.}$$

By similar estimates that lead to (3.15), we can obtain

$$E \sup_{0 \leq t \leq T} \mathbf{e}(u_t^{mn}, v_t^{mn}) \leq C_2 \left\{ \mathbf{e}(g^{mn}, h^{mn}) + E \int_0^T \{\|f_s^{mn}\|^2 + \operatorname{Tr} Q_s^{mn}\} \, ds \right\}$$
(3.17)

for some constant $C_2 > 0$. Since the right-hand side of (3.17) tends to zero as $m, n \to \infty$, it follows that $\|(u^m; v^m) - (u^n; v^n)\| \to 0$. Hence $\{(u^n; v^n)\}$ is a Cauchy sequence in $X$ and $\lim_{n \to \infty}(u^n; v^n) = (u; v)$ strongly as claimed. Due to this strong convergence, it is easy to show that the limit $(u; v)$ is the unique strong solution with the depicted regularity. Moreover, we can take the limits termwise in (3.14) to obtain the energy equation (3.11). Then the energy inequality (3.12) follows easily. □

Notice that, due to the lack of required smoothness of solutions, the general Itô formula does not hold here. As in the deterministic case, the energy equation and the associated inequalities are the key to proving the existence and regularity results for stochastic hyperbolic equations.

Owing to the dissipation term in (3.11), in contrast to the energy inequality (3.12), it is possible to obtain an exponential estimate for the mean energy. To this end, we introduce a *pseudo energy function*

$$(3.18) \quad \mathbf{e}^\lambda(\phi) := \mathbf{e}^\lambda(u; v) = \|Bu\|^2 + \|v + \lambda u\|^2 \quad \text{for } u \in H^1, v \in H,$$

where $\lambda > 0$ is a parameter. Let $v^\lambda = v + \lambda u$. Then we can write

$$(3.19) \quad \mathbf{e}^\lambda(u; v) = \mathbf{e}(u; v^\lambda) = \mathbf{e}(u; v) + 2\lambda(u, v) + \lambda^2 \|u\|^2.$$

Since $A$ is strongly elliptic and strictly positive, its smallest eigenvalue $\eta_1$ can be characterized as ([25], page 62)

$$(3.20) \quad \eta_1 = \inf_{g \in H^1, g \neq 0} \frac{\|g\|_1^2}{\|g\|^2} > 0.$$

LEMMA 3.2 (Exponential estimate). *Let the conditions for Lemma* 3.1 *be satisfied such that* (3.10) *holds for any $T > 0$. Then if*

$$(3.21) \quad \lambda \leq \lambda_0 := \min\left\{\frac{\alpha}{2}, \frac{\eta_1}{4\alpha}\right\},$$

*there exists $\alpha_1 \in (0, \lambda)$ such that the following inequality holds:*

$$(3.22) \quad E\mathbf{e}^\lambda(\phi_t) \leq \mathbf{e}^\lambda(\phi_0)e^{-\alpha_1 t} + \int_0^t e^{-\alpha_1(t-s)} E\left\{\frac{2}{\alpha_1}\|f_s\|^2 + \operatorname{Tr} Q_s\right\} ds.$$

PROOF. It follows from (3.7) that $(u_t; v_t^\lambda)$ satisfies the perturbed system

$$(3.23) \quad \begin{aligned} du_t &= [v_t^\lambda - \lambda u_t]\, dt, \\ dv_t^\lambda &= [Au_t + \lambda(2\alpha - \lambda)u_t - (2\alpha - \lambda)v_t^\lambda + f_t]\, dt + dM_t, \\ u_0 &= g, \quad v_0 = h, \quad 0 < t < T. \end{aligned}$$

By applying Lemma 3.1 to the above system and noting $\mathbf{e}^\lambda(u_t; v_t) = \mathbf{e}(u_t; v_t^\lambda)$, the pseudo energy function (3.18) satisfies

$$
\begin{aligned}
(3.24) \quad d\mathbf{e}^\lambda(u_t; v_t) &= 2[\lambda(2\alpha - \lambda)(u_t, v_t^\lambda) \\
&\quad - \lambda\|u_t\|_1^2 - (2\alpha - \lambda)\|v_t^\lambda\|^2 + (f_t, v_t^\lambda) + \tfrac{1}{2}\operatorname{Tr} Q_t]\, dt \\
&\quad + 2(v_t^\lambda, dM_t),
\end{aligned}
$$

with $\mathbf{e}^\lambda(u_0; v_0) = \mathbf{e}(u_0; v_0^\lambda)$. Now, in view of (3.20), we have, by using some simple inequalities,

$$
\begin{aligned}
\lambda(2\alpha &- \lambda)(u, v^\lambda) - \lambda\|u\|_1^2 - (2\alpha - \lambda)\|v^\lambda\|^2 \\
&\leq \lambda(2\alpha - \lambda)\frac{\|u\|_1}{\sqrt{\eta_1}}\|v^\lambda\| - \lambda\|u\|_1^2 - (2\alpha - \lambda)\|v^\lambda\|^2 \\
&\leq \lambda(2\alpha - \lambda)\left[\lambda\frac{\|u\|_1^2}{\eta_1} + \frac{1}{4\lambda}\|v^\lambda\|^2\right] - \lambda\|u\|_1^2 - (2\alpha - \lambda)\|v^\lambda\|^2 \\
&\leq -\lambda\left(1 - \frac{2\alpha\lambda}{\eta_1}\right)\|u\|_1^2 - \frac{3\alpha}{4}\|v^\lambda\|^2 \leq -\frac{\lambda}{2}\|u\|_1^2 - \frac{3\lambda}{4}\|v^\lambda\|^2.
\end{aligned}
$$

The above result together with the fact that

$$
(f_t, v^\lambda) \leq \frac{\lambda}{4}\|v^\lambda\|^2 + \frac{1}{\lambda}\|f_t\|^2
$$

imply that

$$
\begin{aligned}
(3.25) \quad &\lambda(2\alpha - \lambda)(u, v^\lambda) - \lambda\|u\|_1^2 - (2\alpha - \lambda)\|v^\lambda\|^2 + (f_t, v^\lambda) \\
&\leq -\frac{\lambda}{2}(\|u\|_1^2 + \|v^\lambda\|^2) + \frac{1}{\lambda}\|f_t\|^2.
\end{aligned}
$$

In view of (3.25), equation (3.24) yields

$$
(3.26) \quad d\mathbf{e}^\lambda(u_t; v_t) \leq -\lambda \mathbf{e}^\lambda(u_t; v_t)\, dt + \left[\frac{2}{\lambda}\|f_t\|^2 + \operatorname{Tr} Q_t\right] dt + 2(v_t^\lambda, dM_t),
$$

which can be integrated to get the desired inequality (3.22), after taking the expectation, with any $\alpha_1 < \lambda$. $\square$

It is easy to show that the energy norms induced by $\mathbf{e}$ and $\mathbf{e}^\lambda$ are equivalent. In fact, the following lemma holds.

LEMMA 3.3. *For any $\lambda \in (0, \mu_1)$, the inequality*

$$
(3.27) \quad \left(\frac{\mu_1 - \lambda}{\mu_1 + \lambda}\right)\mathbf{e}(u; v) \leq \mathbf{e}^\lambda(u; v) \leq \left(\frac{\mu_1 + \lambda}{\mu_1 - \lambda}\right)\mathbf{e}(u; v)
$$

holds, where $\mu_1 = (\sqrt{4\eta_1 + \lambda^2})$. Moreover, we have

$$(3.28) \quad E\mathbf{e}(\phi_t) \leq K(\lambda)\left\{\mathbf{e}(\phi_0)e^{-\alpha_1 t} + \int_0^t e^{-\alpha_1(t-s)}E\left(\frac{2}{\alpha_1}\|f_s\|^2 + \operatorname{Tr} Q_s\right)ds\right\},$$

where $K(\lambda) = (\mu_1 + \lambda)/(\mu_1 - \lambda)$.

PROOF. By definition (3.19),
$$\mathbf{e}^\lambda(u;v) = \mathbf{e}(u;v) + 2\lambda(u,v) + \lambda^2\|u\|^2.$$

It follows that, for any $\beta > 0$,

$$\mathbf{e}^\lambda(u;v) \leq \mathbf{e}(u;v) + (1+\beta)\left(\frac{\lambda^2}{\eta_1}\right)\|u\|_1^2 + \frac{1}{\beta}\|v\|^2$$

$$\leq \left(1 + \frac{1}{\beta}\right)\mathbf{e}(u;v) = \left(\frac{\mu_1 + \lambda}{\mu_1 - \lambda}\right)\mathbf{e}(u;v)$$

by choosing $\beta = \frac{1}{2}\{\sqrt{(4\eta_1/\lambda^2) + 1} - 1\}$.

On the other hand, for any $\gamma > 0$,

$$\mathbf{e}^\lambda(u;v) \geq \mathbf{e}(u;v) - \left\{(\gamma - 1)\left(\frac{\lambda^2}{\eta_1}\right)\|u\|_1^2 + \frac{1}{\gamma}\|v\|^2\right\}$$

$$\geq \left(\frac{\mu_1 - \lambda}{\mu_1 + \lambda}\right)\mathbf{e}(u;v)$$

by taking $\gamma = \frac{1}{2}\{\sqrt{1 + 4\eta_1/\lambda^2} + 1\}$.

Therefore, we have verified (3.27), and the result (3.28) is now an direct consequence of (3.22) and (3.27). □

**4. Semilinear stochastic hyperbolic equations.** Let us consider the initial boundary value problem for the hyperbolic equation

$$(4.1) \quad \begin{aligned} \partial_t^2 u(x,t) &= [A(x,D) - 2\alpha\,\partial_t]u(x,t) \\ &\quad + f(u, Du, x, t) + \sigma(u, Du, x, t)\,\partial_t W(x,t), \quad t > 0, \\ u(x,0) &= u_0(x), \qquad \partial_t u(x,0) = v_0(x), \qquad x \in \mathcal{D} \subset \mathbf{R}^d, \\ u(\cdot, t)|_{\partial\mathcal{D}} &= 0, \end{aligned}$$

where, in contrast to the linear problem (3.1), $f(s, y, x, t)$ and $\sigma(s, y, x, t)$ for $x \in \mathcal{D}$, $t > 0$, $s \in \mathbf{R}$ and $y \in \mathbf{R}^d$ are continuous functions, and $W_t = W(\cdot, t)$ is a continuous Wiener random field with covariance operator $R$ with kernel $r(x, y)$ for $x, y \in \mathcal{D}$.

Similarly we rewrite the linear case as a system of Itô equations in $\mathcal{H}^*$:

$$(4.2) \quad \begin{aligned} du_t &= v_t\,dt, \\ dv_t &= [Au_t - 2\alpha v_t + F_t(u_t)]\,dt + dM_t(u) \end{aligned}$$

or

$$u_t = u_0 + \int_0^t u_s \, ds,$$

(4.3)

$$v_t = v_0 + \int_0^t [Au_s - 2\alpha v_s + F_s(u_s)] \, ds + M_t(u),$$

where we set $F_t(u) := f(u, Du, \cdot, t)$,

(4.4) $$M_t(u) = \int_0^t \Sigma_s(u_s) \, dW_s$$

and $\Sigma_t(\cdot) : H^1 \to H$ is defined by $\Sigma_t(u)(x) := \sigma[u(x), Du(x), x, t]$ for any $u \in H^1, x \in \mathcal{D}$.

We are interested in the large-time solutions of (4.1) when the nonlinear terms allow polynomial growth and are locally Lipschitz continuous. For the existence of solutions, we shall impose a set of sufficient conditions. In what follows, for $r, s \in \mathbf{R}$, let $b(r)$ and $k(r,s)$ be real-valued functions which are positive, locally bounded and monotonically increasing in each variable. Let us introduce a positive function $\Theta(\cdot; \cdot) : \mathcal{H} \to \mathbf{R}^+$ which is continuous and locally bounded, and let

(4.5) $$\mathbf{e}(u; v) \leq \Theta(u; v) \leq \mathbf{e}(u; v) + c\|u\|_1^k$$

for any $(u; v) \in \mathcal{H}$ and some constants $c \geq 0$ and $k \geq 2$. As a shorthand notation, we set

$$\|\Sigma_t(u)\|_R^2 = \operatorname{Tr} Q_t(u) = \int_\mathcal{D} r(x,x) [\Sigma_t(u)(x)]^2 \, dx$$

and impose the following conditions, which will be referred to later as Conditions A:

(A1) $A : H^2 \cap H_0^1 \to H$ is an elliptic operator as given in (3.6).
(A2) $F_t(\cdot) : H^1 \to H$ and $\Sigma_t(\cdot) : H^1 \to H$ are continuous in $t \geq 0$. There exist functions $b(r)$ and $k(r,s)$ as indicated above such that, for any $t \geq 0, u \in H^1$,

$$\|F_t(u)\|^2 + \tfrac{1}{2}\|\Sigma_t(u)\|_R^2 \leq b(\|u\|_1) + q(t)$$

for some locally bounded function $q \in L^1(\mathbf{R}^+)$.
(A3) In addition,

$$\|F_t(u) - F_t(u')\|^2 + \tfrac{1}{2}\|\Sigma_t(u) - \Sigma_t(u')\|_R^2 \leq k(\|u\|_1, \|u'\|_1)\|u - u'\|_1^2,$$

for any $u, u' \in H_1, t \geq 0$.

(A4) There exists a positive function $\Theta$ depicted as above and constants $c_i > 0$, $i = 1, 2, 3$, and $\kappa < 1$ such that

$$\int_0^t \{2(F_s(u_s), v_s) + \|\Sigma_s(u_s)\|_R^2\} \, ds$$
$$\leq c_1 + c_2 \int_0^t \Theta(u_s; v_s) \, ds - c_3 \Theta(u_t; v_t) + \kappa \mathbf{e}(u_t; v_t)$$

for any $u. \in \mathcal{C}(\mathcal{R}^+; H^1) \cap \mathcal{C}^1(\mathcal{R}^+; H)$ with $v_t = \partial_t u_t$.

THEOREM 4.1. *Let Conditions A hold true. Then, for $u_0 = g \in H_1$ and $v_0 = h \in H$, the problem (4.1) or the system (4.2) has a unique continuous solution $u. \in \mathcal{C}([0,T]; H_1)$ with $\partial_t u. \in \mathcal{C}([0,T]; H)$ for any $T > 0$. Moreover, the following energy equation holds:*

$$\begin{aligned}
\mathbf{e}(u_t, v_t) &= \mathbf{e}(u_0, v_0) + 2 \int_0^t [(v_s, F_s(u_s)) - 2\alpha \|v_s\|^2] \, ds \\
&\quad + 2 \int_0^t (v_s, \Sigma_s(u_s) \, dW_s) + \int_0^t \|\Sigma_s(u_s)\|_R^2 \, ds \qquad a.s.
\end{aligned}$$
(4.6)

PROOF. The existence proof is similar to that of Theorems 4.1 and 4.2 in [3], and will only be sketched in steps as follows:

*Step* 1. *$H_1$-Lipschitz truncation.* For $N \geq 1$, let $\eta_N(\cdot) : \mathbf{R}^+ = [0, \infty) \to \mathbf{R}^+$ be a $C_0^\infty$ function such that

(4.7) $$\eta_N(s) = \begin{cases} 1, & \text{for } 0 \leq s \leq N/2, \\ 0, & \text{for } s > N, \end{cases}$$

and $0 \leq \eta_N(s) \leq 1$ for $N/2 < s \leq N$. For $(u; v) \in \mathcal{H}$, define $S_N u = \eta_N(\|u\|_1) u$, $F_t^N(u) = \eta_N(\|u\|_1) F_t(S_N u)$ and $\Sigma_t^N(u) = \eta_N(\|u\|_1) \Sigma_t(S_N u)$. Instead of (4.2), consider the truncated system

(4.8) $$\begin{aligned}
du_t &= v_t \, dt, \\
dv_t &= A u_t \, dt + F_t^N(u_t) \, dt + \Sigma_t^N(u) \, dW_t, \\
u_0 &= g, \qquad v_0 = h.
\end{aligned}$$

*Step* 2. *Local solutions.* By conditions (A2) and (A3), it can be shown that

(4.9) $$\|F_t^N(u)\|^2 = \eta_N(\|u\|_1) F_t(S_N u) \leq \alpha_1(N)$$

and

(4.10) $$\|F_t^N(u) - F_t^N(u')\|^2 \leq \alpha_2(N) \|(u - u')\|_1$$

for any $u, u' \in H_1$, $v, v' \in H$ and for some positive constants $\alpha_1, \alpha_2$ depending on $N$. Similarly, we can deduce that

$$\|\Sigma_t^N(u)\|_R^2 \leq \alpha_3(N) \tag{4.11}$$

and

$$\|\Sigma_t^N(u) - \Sigma_t^N(u')\|_R^2 \leq \alpha_4(N)\|J(u - u')\|^2 \tag{4.12}$$

for any $u, u' \in H_1$, $v, v' \in H$, where $\alpha_3, \alpha_4$ are some positive constants depending on $N$.

Therefore the truncated system (4.8) satisfies the usual linear growth and the global Lipschitz condition. By invoking a standard existence theorem ([6], Theorem 7.4), equation (4.3) has a unique solution $\phi^N = (u^N; v^N) \in L^2(\Omega; C([0, T]; H_0^1 \times H))$.

Introduce a stopping time $\tau_N$ defined by

$$\tau_N = \inf\{t > 0 : \|u_t^N\|_1 > N/2\}.$$

Then, for $t < \tau_N$, $u_t = u_t^N$ is the solution of (4.1) with $\partial_t u_t = v_t^N$. As $\tau_N$ is increasing in $N$, let $\tau_\infty = \lim_{N \to \infty} \tau_N$. Define $u_t$ for $t < \tau_\infty \wedge T$ by $u_t = u_t^N$ if $t < \tau_N < T$. Then $u_t$ is the unique local continuous solution.

*Step* 3. *Global solutions.* Assume condition (A4) is also satisfied. By taking the expectation, the energy equation reads

$$E\mathbf{e}(u_{t \wedge \tau_N}; v_{t \wedge \tau_N}) = \mathbf{e}(u_0; v_0) + 2E\int_0^{t \wedge \tau_N}(v_s, F_s(u_s))\,ds$$
$$+ 2E\int_0^{t \wedge \tau_N}(v_s, \Sigma_s(u_s)\,dW_s) + E\int_0^{t \wedge \tau_N}\|\Sigma_s(u_s)\|_R^2\,ds.$$

Letting $\rho_N(t) = E\mathbf{e}(u_{t \wedge \tau_N}; v_{t \wedge \tau_N})$ and invoking condition (A4), the above yields

$$\rho_N(t) \leq [\mathbf{e}(u_0; v_0) + c_1] + c_2\int_0^t \rho_N(s)\,ds + \kappa\rho_N(t). \tag{4.13}$$

Since $\kappa < 1$, there exists $c_3 > 0$ such that

$$\rho_N(T) \leq c_3[\mathbf{e}(u_0; v_0) + c_1]e^{c_2 T} = C_T.$$

On the other hand, we have

$$\rho_N(T) = E\mathbf{e}_\lambda(u_{T \wedge \tau_N}) \geq E\{\mathcal{I}(\tau_N \leq T)\mathbf{e}_\lambda(u_{T \wedge \tau_N})\}$$
$$\geq CE\{\|u_{T \wedge \tau_N}\|_1^2 \mathcal{I}(\tau_N \leq T)\} \geq C\left(\frac{N}{2}\right)^2 P\{\tau_N \leq T\},$$

where $\mathcal{I}$ is the indicator function and $C > 0$ is a constant. The above inequality gives

$$P\{\tau_N \leq T\} \leq 4\rho_N(T)/CN^2 \leq 4C_T/CN^2.$$

Since the series $\sum_{N=1}^{\infty} P\{\tau_N \leq T\}$ converges, by the Borel–Cantelli lemma, we can conclude that

$$P\{\tau_\infty \leq T\} = 0$$

or $\tau_\infty > T$ a.s. for any $T > 0$. Now we let $u_t^N = u_{t \wedge \tau_N}$ and denote its limit $\lim_{N \to \infty} u_t^N$ still by $u_t$. Then $u_t$ is the global solution as claimed. The energy equation (4.6) can be verified by taking the limits termwise, as $N \to \infty$, in the energy equation for the $N$th truncated system (4.8). $\square$

REMARK. In the above theorem, for simplicity, we assumed that $W(x,t)$ is a scalar Wiener random field. Under an obvious modification, Theorem 4.1 and the subsequent theorems still hold true when $W = (W^{(1)}, \ldots, W^{(k)})$ is a $k$-vector-valued Wiener random field and $\sigma = (\sigma_1, \ldots, \sigma_k)$ is another $k$-vector-valued, predictable random field such that the product $\sigma(\cdot)W(\cdot) = \sum_{j=1}^{k} \sigma_j W^{(j)}(\cdot)$ is interpreted as a dot product.

**5. Bounded solutions.** In view of (4.2), we rewrite the hyperbolic system (4.2) as

(5.1)
$$\begin{aligned} du_t &= v_t \, dt, \\ dv_t &= [Au_t - 2\alpha v_t + F_t(u_t)] \, dt + \Sigma_t(u_t) \, dW_t, \qquad t > 0, \end{aligned}$$

with a given initial state $(u_0; v_0)$ which is an $\mathcal{F}_0$ random vector in $\mathcal{H}$. For the existence of bounded solutions, we shall impose Conditions B as follows:

(B1) There exist $\Phi \in \mathcal{C}^1(H^1; \mathbf{R}^+)$ with Fréchet derivative $\Phi' \in \mathcal{C}(H^1; H)$ and $p_\cdot \in \mathcal{C}(\mathbf{R}^+ \times H^1; H)$ such that $F_t(u) = -\frac{1}{2}\Phi'(u) + p_t(u)$ for any $u \in H^1$ and

$$c_1 \leq \Phi(u) \leq c_2(1 + \|u\|_1^k)$$

for some constants $c_1$ and $c_2 > 0, k \geq 2$.

(B2) There exist constants $\beta_i \geq 0$ and $\gamma_i, \delta_1 \in \mathbf{R}$ with $i = 1, 2, 3$, and essentially bounded functions $\theta$ and $\rho$ which are locally integrable such that

$$\begin{aligned} (\Phi'(u), u) &\geq \beta_1 \Phi(u) - \gamma_1 \|u\|_1^2 - \delta_1, \\ \|p_t(u)\|^2 &\leq \beta_2 \Phi(u) + \gamma_2 \|u\|_1^2 + \theta(t) \end{aligned}$$

and

$$\|\Sigma_t(u)\|_R^2 \leq \beta_3 \Phi(u) + \gamma_3 \|u\|_1^2 + \rho(t)$$

for any $u \in H_1$ and $t > 0$.

(B3) The above constants satisfy

$$(\beta_1 - \tfrac{1}{2})\lambda^2 - \beta_3\lambda - \beta_2 \geq 0,$$
$$(\gamma_1 - \tfrac{1}{2})\lambda^2 + \gamma_3\lambda + \gamma_2 \leq 0.$$

THEOREM 5.1 (Bounded in mean-square). *Suppose that Conditions A and B hold true. Given $u_0 \in H^1$ and $v_0 \in H$ being $\mathcal{F}_0$ random variables such that*

$$E\{\mathbf{e}(u_0; v_0) + \Phi(u_0)\} < \infty,$$

*then the solution of the problem* (5.1) *is bounded in mean-square. Moreover, there exist positive constants $K_1$ and $\alpha_2 > 0$ such that*

(5.2)
$$\begin{aligned} E\{\mathbf{e}(u_t; v_t) + \Phi(u_t)\} \\ \leq K_1 \bigg\{ & E[\mathbf{e}(u_0; v_0) + \Phi(u_0; v_0)]e^{-\alpha_2 t} \\ & + \int_0^t e^{-\alpha_2(t-s)}\left[\frac{1}{\lambda}\theta(s) + \rho(s)\right]ds\bigg\} + 2|\delta_1| \qquad \forall t > 0. \end{aligned}$$

PROOF. By applying Lemma 3.2 to (5.1), as in (3.24), we obtain the perturbed energy equation

(5.3)
$$\begin{aligned} d\mathbf{e}^\lambda(u_t; v_t) = 2[&\lambda(2\alpha - \lambda)(u_t, v_t^\lambda) - \lambda\|u_t\|_1^2 - (2\alpha - \lambda)\|v_t^\lambda\|^2 \\ & + (F_t(u_t), v_t^\lambda) + \tfrac{1}{2}\|\Sigma_t(u_t)\|_R^2]\,dt \\ & + 2\,dM_t^\lambda(u; v), \end{aligned}$$

where we set $dM_t^\lambda(u; v) = (v_t^\lambda, \Sigma(u_t)\,dW_t)$. As in (3.27) in Lemma 3.2, for $\lambda < \{\alpha \wedge \eta_1/4\alpha\}$, the above yields

(5.4)
$$\begin{aligned} d\mathbf{e}^\lambda(u_t; v_t) \leq & -\lambda(\|u_t\|_1^2 + \tfrac{3}{2}\|v_t^\lambda\|^2)\,dt \\ & + [2(F_t(u_t), v_t^\lambda) + \|\Sigma_t(u_t)\|_R^2]\,dt + 2\,dM_t^\lambda(u; v). \end{aligned}$$

By assumptions,

$$\begin{aligned} (F_t(u_t), v_t^\lambda) &= -\frac{1}{2}(\Phi'(u_t), v_t) - \frac{1}{2}\lambda(\Phi'(u_t), u_t) + (p_t(u_t), v_t^\lambda) \\ &\leq -\frac{1}{2}\frac{d}{dt}\Phi(u_t) - \frac{1}{2}\lambda(\Phi'(u_t), u_t) + \frac{\lambda}{2}\|v_t^\lambda\|^2 + \frac{1}{2\lambda}\|p_t(u_t)\|^2, \end{aligned}$$

which, in view of condition (B2), implies that

(5.5)
$$\begin{aligned} (F_t(u_t), v_t^\lambda) \leq & -\frac{1}{2}\frac{d}{dt}\Phi(u_t) - \frac{1}{2}\left(\beta_1\lambda - \frac{\beta_2}{\lambda}\right)\Phi(u_t) \\ & + \frac{\lambda}{2}\|v_t^\lambda\|^2 + \frac{1}{2}\left(\gamma_1\lambda + \frac{\gamma_2}{\lambda}\right)\|u_t\|_1^2 + \frac{1}{2\lambda}\theta(t) + \frac{\lambda}{2}\delta_1. \end{aligned}$$

Define a *superenergy function* $J: \mathcal{H} \to \mathcal{R}^+$ by

(5.6) $$J(u;v) = \mathbf{e}(u;v) + \Phi(u),$$

with $J^\lambda = \mathbf{e}^\lambda + \Phi$. By applying (5.4), (5.5) and condition (B2) to (5.3), we obtain

$$dJ^\lambda(u_t;v_t) \leq -\lambda \mathbf{e}^\lambda(u_t;v_t)\,dt - \left(\beta_1 \lambda - \frac{\beta_2}{\lambda} - \beta_3\right)\Phi(u_t)\,dt$$
$$+ \left\{\frac{\lambda}{2}\|v_t^\lambda\|^2 + \left(\gamma_1 \lambda + \frac{\gamma_2}{\lambda} + \gamma_3\right)\|u_t\|_1^2 + \frac{1}{\lambda}\theta(t) + \rho(t) + \lambda \delta_1\right\} dt$$
$$+ 2\,dM_t^\lambda(u;v).$$

By invoking condition (B3), the above inequality gives

(5.7) $$dJ^\lambda(u_t;v_t) \leq -\frac{\lambda}{2} J^\lambda(u_t;v_t)\,dt + \left\{\frac{1}{\lambda}\theta(t) + \rho(t) + \lambda\delta_1\right\} dt$$
$$+ 2\,dM_t^\lambda(u;v),$$

which implies that

(5.8) $$EJ^\lambda(u_t;v_t) \leq EJ^\lambda(u_0;v_0)e^{-\lambda t/2}$$
$$+ \int_0^t e^{-\lambda(t-s)/2}\left[\frac{1}{\lambda}\theta(s) + \rho(s)\right] ds + 2|\delta_1| < \infty$$

for all $t > 0$. Since, by assumption, $\theta$ and $\rho \in L^1_{\text{loc}}(\mathbf{R}^+)$ are essentially bounded, we have

$$EJ^\lambda(u_t;v_t) = E\{\mathbf{e}^\lambda(u_t;v_t) + \Phi(u_t;v_t)\} < \infty.$$

Now, by invoking Lemma 3.3, $J(u;v) \leq CJ^\lambda(u;v)$ for some $C > 0$. Therefore, the result (5.2) holds with some constant $K_1 > 0$ and $\alpha_2 = \frac{\lambda}{2}$. $\square$

In fact, under somewhat stronger assumptions, it is possible to show that the solution of (5.1) is ultimately bounded in mean-square, that is,

$$\limsup_{t \to \infty} E\{\|u_t\|_1^2 + \|v_t\|^2\} < \infty.$$

THEOREM 5.2. *Assume that Conditions* A *and* B *hold true with* $\delta_1 = 0, \theta$ *and* $\rho \in L^1(\mathbf{R}^+)$. *Then the solution* $\phi_t = (u_t;v_t)$ *is ultimately bounded in mean-square such that*

(5.9) $$E \sup_{0 \leq t \leq T}\{\mathbf{e}(u_t;v_t) + \Phi(u_t)\}$$
$$\leq K_2 E\{\mathbf{e}(u_0;v_0) + \Phi(u_0)\} + K_3 \int_0^T [\theta(s) + \rho(s)]\,ds$$

*for some positive constants* $K_2$ *and* $K_3$.

PROOF. In view of (5.7), it is clear that

$$J^\lambda(u_t; v_t) + \frac{\lambda}{2}\int_0^t J^\lambda(u_s; v_s)\, ds$$

$$\leq J^\lambda(u_0; v_0) + \int_0^t \left[\frac{1}{\lambda}\theta(s) + \rho(s)\right] ds + 2M_t^\lambda(u; v).$$

Hence

(5.10)
$$EJ^\lambda(u_t; v_t) + \frac{\lambda}{2}\int_0^t EJ^\lambda(u_s; v_s)\, ds$$
$$\leq EJ^\lambda(u_0; v_0) + \int_0^t \left[\frac{1}{\lambda}\theta(s) + \rho(s)\right] ds$$

and

(5.11)
$$E\sup_{0\leq t\leq T} J^\lambda(u_t; v_t) \leq EJ^\lambda(u_0; v_0) + \int_0^t \left[\frac{1}{\lambda}\theta(s) + \rho(s)\right] ds$$
$$+ 2E\sup_{0\leq t\leq T}|M_t^\lambda(u; v)|.$$

By means of the B–D–G inequality for a submartingale, we can deduce that

(5.12)
$$E\sup_{0\leq t\leq T}|M_t^\lambda(u; v)| = E\sup_{0\leq t\leq T}\left|\int_0^t (v_s^\lambda, \Sigma_s(u_s)\, dW_s)\right|$$
$$\leq 3E\left\{\int_0^T (R\Sigma_s(u_s)v_s^\lambda, \Sigma_s(u_s)v_s^\lambda)\, ds\right\}^{1/2}$$
$$\leq 3E\left\{\sup_{0\leq t\leq T}\|v_t^\lambda\|\right\}\left\{\int_0^T \|\Sigma_s(u_s)\|_R^2\, ds\right\}^{1/2}$$
$$\leq \tfrac{1}{4}E\sup_{0\leq t\leq T}\|v_t^\lambda\|^2 + 9E\int_0^T \|\Sigma_s(u_s)\|_R^2\, ds.$$

Results (5.8) and (5.12) and condition (B2) imply that

$$E\sup_{0\leq t\leq T} J^\lambda(u_t; v_t) \leq EJ^\lambda(u_0; v_0) + \frac{1}{2}E\sup_{0\leq t\leq T}\|v_t^\lambda\|^2$$
$$+ 18E\int_0^T \|\Sigma_s(u_s)\|_R^2\, ds + \int_0^t \left[\frac{1}{\lambda}\theta(s) + \rho(s)\right] ds$$
$$\leq EJ^\lambda(u_0; v_0) + \frac{1}{2}E\sup_{0\leq t\leq T}\|v_t^\lambda\|^2$$
$$+ 18E\int_0^T [\beta_3\Phi(u_s) + \gamma_3\|u_s\|_1^2]\, ds$$

$$+ \int_0^t \left[\frac{1}{\lambda}\theta(s) + 19\rho(s)\right] ds.$$

Therefore, there exist positive constants $c_i, i = 1, 2, 3$, such that

$$E \sup_{0 \leq t \leq T} J^\lambda(u_t; v_t) \leq c_1 E J^\lambda(u_0; v_0) + c_2 \int_0^T E \sup_{0 \leq \tau \leq s} J^\lambda(u_\tau; v_\tau) \, ds$$

$$+ c_3 \int_0^T [\theta(s) + \rho(s)] \, ds.$$

From this together with the bound (5.10) and the Gronwall lemma, we can infer that there exists a pair of positive constants $k_2, k_3$ such that

$$E \sup_{0 \leq t \leq T} J^\lambda(u_t; v_t) \leq k_2 E J^\lambda(u_0; v_0) + k_3 \int_0^T [\theta(s) + \rho(s)] \, ds,$$

which, by Lemma 3.3, leads to the desired inequality (5.9). $\square$

**6. Asymptotic stability of solutions.** Suppose that the hyperbolic system (5.1) has an equilibrium solution $u = \hat{u} \in \mathcal{D}(A)$ with $v = 0$. By a translation via $\hat{u}$, without loss of generality, we may assume that $(\hat{u}; \hat{v}) \equiv (0; 0)$ is an equilibrium solution. We are interested in the asymptotic stability of the null solution in the following sense.

DEFINITIONS.

1. The null solution $\phi = (u; v) \equiv (0; 0)$ of (5.1) is said to be *asymptotically stable in mean-square* in $\mathcal{H}$ if $\exists \delta > 0$ such that, for $\|\phi_0\| < \delta$,

$$\lim_{t \to \infty} E\|\phi_t\|^2 = 0,$$

and it is *exponentially stable in mean-square* if there exist positive constants $K(\delta)$ and $\nu$ such that

$$E\|\phi_t\|^2 \leq K(\delta) e^{-\nu t} \qquad \forall t > 0,$$

where $\|\phi\|^2 = \|u\|_1^2 + \|v\|^2$.

2. The null solution is said to be a.s. (almost surely) *asymptotically stable* if

$$P\left\{\lim_{t \to \infty} \|\phi_t\| = 0\right\} = 1,$$

and it is *a.s. exponentially stable* if there exist positive constants $K_2(\delta), \nu_2$ and a random time $T(\omega) > 0$ such that

$$\|\phi_t\| \leq K_2(\delta) e^{-\nu_2 t} \qquad \forall t > T, \text{ a.s.}$$

REMARK. In view of the above definitions, it is clear that the exponential stability implies the asymptotic stability.

To proceed we assume that $F_t(0) = 0$ and $\Sigma_t(0) = 0$ for any $t > 0$ so that $\phi = (u; v) \equiv (0; 0)$ is an equilibrium solution of equation (5.1). In the stability analysis [11], it is often assumed that the global solution exists in the first place. Hence we suppose, under suitable conditions such as Conditions A, that the equation has a unique global solution.

THEOREM 6.1 (Stability in mean-square). *Suppose that Conditions B hold true with the following provisions:*

(1) $\Phi(0) = 0$ *and* $\Phi(u) > 0$ *if* $u \neq 0$.
(2) *In condition* (B2), $\delta_1 = 0$ *and there exists* $\alpha_0 > 0$ *such that*

$$\int_0^\infty e^{\alpha_0 t} \left[\frac{1}{\lambda}\theta(t) + \rho(t)\right] dt = C < \infty.$$

*Then the null solution of equation* (5.1) *is exponentially stable in mean-square. Moreover, if* $\phi_0 = (u_0; v_0)$ *is an* $\mathcal{F}_0$-*measurable random variable in* $\mathcal{H}$ *satisfying*

$$E\{\mathbf{e}(u_0; v_0) + \Phi(u_0; v_0)\} < \infty,$$

*then the inequality*

(6.1)     $E\{\mathbf{e}(u_t; v_t) + \Phi(u_t)\} \leq E\{\mathbf{e}(u_0; v_0) + \Phi(u_0) + C\}e^{-\alpha_3 t}$

*holds for any* $t > 0$, *where* $\alpha_3 = (\alpha_0 \wedge \alpha_2)$.

The theorem follows immediately from (5.2) in Theorem 5.1 and the simple fact that $\int_0^\infty e^{-\alpha_2(t-s)}\theta(s)\,ds \leq e^{-\alpha_0 t}\int_0^\infty e^{\alpha_0 s}\theta(s)\,ds$. In fact it is possible to show that the null solution is a.s. exponentially stable. Before stating the next theorem, we need a lemma which is a simple consequence of Theorem 5.2.

LEMMA 6.2. *Under the conditions for Theorem* 5.1, *the solution* $\phi_t = (u_t; v_t)$ *is ultimately bounded in mean-square such that*

(6.2)
$$E \sup_{0 \leq t \leq T} \{\mathbf{e}^\lambda(u_t; v_t) + \Phi(u_t)\}$$
$$\leq K_1 E\{\mathbf{e}^\lambda(u_0; v_0) + \Phi(u_0)\} + K_2 \int_0^T [\theta(t) + \rho(t)]\,dt$$

*for some constants* $K_1, K_2 > 0$.

With the aid of Theorem 6.1 and Lemma 6.2, we can prove the following theorem.

THEOREM 6.3 (Almost sure stability). *Assume that all of the conditions for Theorem* 5.1 *hold true with* $\theta = \rho \equiv 0$. *Then the null solution of* (5.1) *is exponentially stable almost surely. Moreover, there exist positive constants* $C, \nu$ *and a random variable* $T(\omega) > 0$ *such that*

(6.3) $\quad \mathbf{e}(u_t; v_t) + \Phi(u_t) \leq C\{\mathbf{e}(u_0; v_0) + \Phi(u_0)\}e^{-\nu t} \qquad a.s.$

*for any* $t > T$.

PROOF. Owing to Lemma 3.3, instead of (6.3), it suffices to show that

(6.4) $\quad \{\mathbf{e}^\lambda(u_t; v_t) + \Phi(u_t)\} \leq C_0 E\{\mathbf{e}^\lambda(u_0; v_0) + \Phi(u_0)\}e^{-\nu t}, \qquad t > T, \text{ a.s.}$

for some constant $C_0 > 0$ and for $\lambda$ satisfying (3.21). To this end, decompose $\mathbf{R}^+$ as $\mathbf{R}^+ = \bigcup_{n=0}^\infty [n, n+1]$ and consider the solution $(u_t; v_t)$ for $n \leq t < n+1$. Following the steps leading to (5.8) in the proof of Theorem 5.1, it can be shown that

(6.5) $\quad E \sup_{n \leq t \leq n+1} J^\lambda(u_t; v_t) \leq E J^\lambda(u_n; v_n) + 2E \sup_{0 \leq t \leq T} |M_t^\lambda(u)|,$

where we recall that $J^\lambda(u; v) = \mathbf{e}^\lambda(u; v) + \Phi(u)$. As in (5.12), we have

$$E \sup_{n \leq t \leq n+1} |M_t^\lambda(u; v)|$$

$$\leq 3E \left\{ \sup_{n \leq t \leq n+1} \|v_t^\lambda\| \right\} \left\{ \int_n^{n+1} \|\Sigma_s(u_s)\|_R^2 ds \right\}^{1/2}$$

$$\leq 3 \left\{ E \sup_{n \leq t \leq n+1} \|v_t^\lambda\|^2 \right\}^{1/2} \left\{ E \int_n^{n+1} \|\Sigma_s(u_s)\|_R^2 ds \right\}^{1/2}$$

$$\leq 3 \left\{ E \sup_{n \leq t \leq n+1} \|v_t^\lambda\|^2 \right\}^{1/2} \left\{ (\beta_3 \vee \gamma_3) \int_n^{n+1} E J^\lambda(u_s; v_s) ds \right\}^{1/2}.$$

By making use of (5.11), the above gives rise to the upper bound

$$E \sup_{n \leq t \leq n+1} |M_t^\lambda(u; v)|$$

(6.6) $\quad \leq 3\{EJ^\lambda(u_0; v_0)\} \left\{ CK(\beta_3 \vee \gamma_3) \int_n^{n+1} e^{-\lambda s/2} ds \right\}^{1/2}$

$$\leq 6\{EJ^\lambda(u_0; v_0)\} \left\{ \frac{CK}{\lambda}(\beta_3 \vee \gamma_3) \right\}^{1/2} e^{-n\lambda/4}.$$

By taking (6.1), (6.5) and (6.6) into account, we get

(6.7) $\quad E \left\{ \sup_{n \leq t \leq n+1} J^\lambda(u_t; v_t) \right\} \leq C_0 \{EJ^\lambda(u_0; v_0)\} e^{-n\lambda/4}$

for constant $C_0 > 0$.

Therefore, by using the Markov inequality and (6.7),

$$P\left\{\sup_{n \leq t \leq n+1} J^\lambda(u_t; v_t) > C_0 E J^\lambda(u_0; v_0) e^{-n\lambda/8}\right\}$$

$$\leq \frac{E\{\sup_{n \leq t \leq n+1} J^\lambda(u_t; v_t)\}}{C_0 E\{J^\lambda(u_0; v_0)\} e^{-n\lambda/8}} \leq e^{-n\lambda/8}.$$

Since $\sum_{n=1}^\infty e^{-n\lambda/8} < \infty$, it follows from the Borel–Cantelli lemma that there exists a random number $N(\omega) > 0$ such that, for $n > N$,

$$\sup_{n \leq t \leq n+1} J^\lambda(u_t; v_t) \leq C_0 \{E J^\lambda(u_0; v_0)\} e^{-n\lambda/8} \qquad \text{a.s.},$$

which, by definition (5.6), implies (6.4) with $\nu = \lambda/8$.  $\square$

**7. Invariant measures.** Let us consider the autonomous version of the system (5.1):

(7.1)
$$du_t = v_t\, dt,$$
$$dv_t = [Au_t - 2\alpha v_t + F(u_t)]\, dt + \Sigma(u_t)\, dW_t, \qquad t > 0,$$

with a given initial state $(u_0; v_0)$, where $F$ and $\Sigma$ do not depend on $t$ explicitly.

Let $\phi_t = (u_t; v_t)$ and rewrite the system (7.1) as an evolution equation in the differential form,

(7.2) $$d\phi_t = \mathcal{A}\phi_t\, dt + \mathcal{F}(\phi_t)\, dt + d\mathcal{M}_t(\phi)$$

with $\phi_0 = (u_0; v_0)$, where we set

$$\phi_t = \begin{bmatrix} u_t \\ v_t \end{bmatrix}, \qquad \mathcal{F}(\phi) = \begin{bmatrix} 0 \\ F(u) \end{bmatrix}, \qquad \mathcal{M}_t(\phi) = \begin{bmatrix} 0 \\ M_t(u) \end{bmatrix}$$

and

$$\mathcal{A} = \begin{bmatrix} 0 & I \\ A & -2\alpha I \end{bmatrix}$$

where $I$ is the identity operator on $H$.

Under Conditions C given below, as an Itô equation in a Hilbert space, the solution $\phi_t$, if it exists, is a Markov diffusion process in $\mathcal{H}$ (see [6], Chapter 9). The transition probability function is given by

$$P_t(\xi; \mathcal{B}) = P\{\phi_t \in \mathcal{B} | \phi_0 = \xi\}, \qquad \xi \in \mathcal{H},\ \mathcal{B} \in \sigma(\mathcal{H}).$$

Suppose there exists an invariant measure $\mu$ on $(\mathcal{H}, \sigma(\mathcal{H}))$, where $\sigma(\mathcal{H})$ denotes the Borel $\sigma$-field of $\mathcal{H}$. Then it satisfies ([7], page 12):

$$\mu(\mathcal{B}) = \int_\mathcal{H} P_t(\xi; \mathcal{B}) \mu(d\xi) \qquad \forall \mathcal{B} \in \sigma(\mathcal{H}).$$

To show the existence of an invariant measure, we shall specialize Conditions A by assuming that the nonlinear terms satisfy a uniform Lipschitz continuity condition. To be precise, assume the following Conditions C:

(C1) Let $F(\cdot)\colon H^1 \to H$ and $\Sigma(\cdot)\colon H^1 \to H$, and let there exist positive constants $b_i, c_i$ for $i = 1, 2$, such that
$$\|F(u)\|^2 \leq b_1 \|u\|_1^2 + c_1$$
and
$$\|\Sigma_t(u)\|_R^2 \leq b_2 \|u\|_1^2 + c_2$$
for any $u \in H^1$.

(C2) There exist positive constants $k_1, k_2$ such that
$$\|F(u) - F(u')\|^2 \leq k_1 \|u - u'\|_1^2$$
and
$$\|\Sigma(u) - \Sigma(u')\|_R^2 \leq k_2 \|u - u'\|_1^2$$
for any $u, u' \in H^1$.

(C3) The constants $b_i$ and $k_i$ satisfy
$$(b_1 + b_2 \lambda) \wedge (k_1 + k_2 \lambda) \leq \frac{\lambda^2}{2}.$$

To show the existence of an invariant measure, we shall follow an approach by Da Prato and Zabczyk ([7], Theorem 6.3.2) for some stochastic dissipative systems. Though not directly applicable to the present problem, it can be adapted to proving the following theorem.

THEOREM 7.1 (Invariant measures). *Suppose that the system* (7.1) *satisfies Conditions* C. *Then there exists a unique invariant measure $\mu$ on $(\mathcal{H}, \sigma(\mathcal{H}))$. Moreover, given any bounded Lipschitz continuous function $G$ on $\mathcal{H}$, there are positive constants $C$ and $\alpha_2$ such that*

$$(7.3) \qquad \left| \int_{\mathcal{H}} G(\eta) P_t(\xi; d\eta) - \int_{\mathcal{H}} G(\eta) \mu(d\eta) \right| \leq C(1 + \|\xi\|) e^{-\alpha_2 t}$$

*for any $t > 0$ and $\xi \in \mathcal{H}$.*

PROOF. To extend the time domain for the system (7.1) to the whole real line **R**, introduce an independent copy $V_t$ of the Wiener process $W_t$ for $t \geq 0$. Define $\hat{W}_t$ by

$$(7.4) \qquad \hat{W}_t = \begin{cases} W_t, & \text{for } t \geq 0, \\ V_{-t}, & \text{for } t \leq 0, \end{cases}$$

and let $\hat{\mathcal{F}}_t = \sigma\{\hat{W}_s : s \leq t\}$ for $t \in \mathbf{R}$. Now, for $t > \tau$, let $\phi_t(\tau;\xi) = (u_t; v_t)(\tau;\xi) = (u_t(\tau;\xi); v_t(\tau;\xi))$ be the solution of the extended system

$$
\begin{aligned}
du_t &= v_t \, dt, \\
dv_t &= [Au_t - 2\alpha v_t + F(u_t)] \, dt + \Sigma(u_t) \, d\hat{W}_t, \quad t > \tau, \\
u_\tau &= \xi_1, \qquad v_\tau = \xi_2,
\end{aligned}
\tag{7.5}
$$

where $\xi = (\xi_1; \xi_2) \in \mathcal{H}$.

Similarly to the derivation of the inequality (5.4) in the proof of Theorem 5.1, it can be shown that, for $\lambda < \{\alpha \wedge \eta_1/4\alpha\}$,

$$
d\mathbf{e}^\lambda(u_t; v_t) \leq -\lambda \mathbf{e}^\lambda(u_t; v_t) \, dt + \left[ \frac{\lambda}{2} \|v_t^\lambda\|^2 + \frac{1}{\lambda}\|F(u_t)\|^2 + \|\Sigma(u_t)\|_R^2 \right] dt
$$
$$
+ 2(v_t^\lambda, \Sigma(u_t) \, d\hat{W}_t),
\tag{7.6}
$$

where $d\hat{M}_t^\lambda(u;v) = (v_t^\lambda, \Sigma(u_t) \, d\hat{W}_t)$. By making use of conditions (C1) and (C2), the above yields

$$
\begin{aligned}
d\mathbf{e}^\lambda(u_t; v_t) &\leq -\lambda \mathbf{e}^\lambda(u_t; v_t) \, dt + \left[ \frac{\lambda}{2}\|v_t^\lambda\|^2 + \left(\frac{b_1}{\lambda} + b_2\right)\|u_t\|_1^2 + \left(\frac{c_1}{\lambda} + c_2\right) \right] dt \\
&\quad + 2 \, d\hat{M}_t^\lambda(u;v) \\
&\leq \left[ -\frac{\lambda}{2}\mathbf{e}^\lambda(u_t; v_t) + \left(\frac{c_1}{\lambda} + c_2\right) \right] dt + 2 \, d\hat{M}_t^\lambda(u;v),
\end{aligned}
$$

so that

$$
E\mathbf{e}^\lambda[\phi_t(s,\xi)] \leq \mathbf{e}^\lambda(\xi) e^{-\lambda(t-s)/2} + \left(\frac{c_1}{\lambda} + c_2\right) \int_s^t e^{-\lambda(t-r)/2} \, dr.
$$

Therefore, there exists a constant $K_1 > 0$ such that

$$
E\mathbf{e}^\lambda[\phi_t(s,\xi)] \leq K_1\{1 + \mathbf{e}^\lambda(\xi)\} \qquad \text{for any } t \geq s.
\tag{7.7}
$$

For $\tau_1 > \tau_2 > 0$, let

$$
u_t^i = u_t(-\tau_i; \xi), \qquad v_t^i = v_t(-\tau_i; \xi) \qquad \text{for} \quad t > -\tau_2,
$$

with $i = 1, 2$ and

$$
\tilde{u}_t = u_t^1 - u_t^2, \qquad \tilde{v}_t = v_t^1 - v_t^2.
$$

Then it follows from (7.5) that we have

$$
\begin{aligned}
d\tilde{u}_t &= \tilde{v}_t \, dt, \\
d\tilde{v}_t &= [A\tilde{v}_t - 2\alpha\tilde{v}_t + \delta F(u_t^1; u_t^2)] \, dt + \delta\Sigma(u_t^1; u_t^2) \, d\hat{W}_t, \quad t > -\tau_2, \\
\tilde{u}_{-\tau_2} &= (u_{-\tau_2}^1 - \xi_1), \qquad \tilde{v}_{-\tau_2} = (v_{-\tau_2}^1 - \xi_2),
\end{aligned}
\tag{7.8}
$$

where

(7.9)
$$\delta F(u^1; u^2) = F(u^1) - F(u^2),$$
$$\delta \Sigma(u^1; u^2) = \Sigma(u^1) - \Sigma(u^2).$$

Let $\tilde{v}^\lambda = \tilde{v} + \lambda \tilde{u}$ and $\tilde{\mathbf{e}}_t^\lambda = \mathbf{e}(\tilde{u}_t; \tilde{v}_t^\lambda)$. As in (7.6), we can obtain the energy inequality

(7.10)
$$d\tilde{\mathbf{e}}_t^\lambda \leq -\lambda \tilde{\mathbf{e}}_t^\lambda \, dt + \left[ \frac{\lambda}{2} \|\tilde{v}_t^\lambda\|^2 + \frac{1}{\lambda} \|\delta F(u_t^1; u_t^2)\|^2 + \|\delta \Sigma(u_t^1; u_t^2)\|_R^2 \right] dt$$
$$+ 2(\tilde{v}_t^\lambda, \delta \Sigma(u_t) \, d\hat{W}_t).$$

In view of (7.9) and conditions (C2) and (C3), equation (7.10) yields

$$d\tilde{\mathbf{e}}_t^\lambda \leq -\lambda \tilde{\mathbf{e}}_t^\lambda \, dt + \left[ \frac{\lambda}{2} \|\tilde{v}_t^\lambda\|^2 + \left( \frac{k_1}{\lambda} + k_2 \right) \|\tilde{u}_t\|_1^2 \right] dt + 2(\tilde{v}_t^\lambda, \delta \Sigma(u_t) \, d\hat{W}_t)$$
$$\leq -\frac{\lambda}{2} \tilde{\mathbf{e}}_t^\lambda \, dt + 2(\tilde{v}_t^\lambda, \delta \Sigma(u_t) \, d\hat{W}_t),$$

which implies that

$$E\tilde{\mathbf{e}}_t^\lambda \leq E\tilde{\mathbf{e}}_{-\tau_2}^\lambda e^{-\lambda(t+\tau_2)/2}.$$

In view of the bound (7.7) and the initial conditions in (7.8), we can show that

$$E\tilde{\mathbf{e}}_{-\tau_2}^\lambda \leq 2K_1\{1 + \mathbf{e}^\lambda(\xi)\},$$

so that

(7.11) $\quad E\tilde{\mathbf{e}}_t^\lambda = E\mathbf{e}^\lambda(u_t^1 - u_t^2; v_t^1 - v_t^2) \leq 2K_1\{1 + \mathbf{e}^\lambda(\xi)\} e^{-\lambda(t+\tau_2)/2}.$

Let

$$\psi_\tau = (u_0(-\tau; \xi); v_0(-\tau; \xi)).$$

By setting $t = 0$ in (7.11), we obtain, for any $\xi \in \mathcal{H}$,

(7.12) $\quad E\mathbf{e}^\lambda(\psi_{\tau_2} - \psi_{\tau_1}) \leq 2K_1\{1 + \mathbf{e}^\lambda(\xi)\} e^{-\lambda \tau_2/2},$

which goes to zero as $\tau_2 \to \infty$. By Lemma 3.3 with $\lambda > 0$, the energy functions $\mathbf{e}^\lambda$ and $\mathbf{e}$ are equivalent so that $\mathbf{e}(\cdot) \leq C\mathbf{e}^\lambda(\cdot)$ for some constant $C > 0$. Since $\mathbf{e}$ defines the energy norm $\|\cdot\|$ on $\mathcal{H}$ by $\|\phi\|^2 = \mathbf{e}(\phi), \phi \in \mathcal{H}$, the set $\{\psi_\tau : \tau \geq 0\}$ is a Cauchy family of random variables in $L^2(\Omega; \mathcal{H})$. Therefore, there exists a unique random variable $\psi_\infty \in L^2(\Omega; \mathcal{H})$ such that, for any $\xi \in \mathcal{H}$,

$$\lim_{\tau \to \infty} E\|\psi_\tau - \psi_\infty\|^2 = 0.$$

Since the random variables $\psi_\tau = \phi_\tau(0;\xi)$ in distribution, $\phi_\tau(0;\xi)$ converges weakly to $\psi_\infty$ as $\tau \to \infty$. It follows that $P_t(\xi;\cdot)$ converges weakly to the probability measure $\mu$ for $\psi_\infty$ as $t \to \infty$, and $\mu$ is the desired invariant measure for the system (7.1).

To verify (7.3), let $G: \mathcal{H} \to \mathbf{R}$ be bounded and Lipschitz continuous on $\mathcal{H}$ such that
$$|G(\xi) - G(\eta)| \leq \gamma \|\xi - \eta\|.$$

Then, for any $t > s > 0$, we have
$$\left| \int_\mathcal{H} G(\eta) P_t(\xi;d\eta) - \int_\mathcal{H} G(\eta) P_s(\xi;d\eta) \right|$$
$$= |EG(\psi_t) - EG(\psi_s)|$$
$$\leq \gamma \{E\|\psi_t - \psi_s\|^2\}^{1/2}$$
$$\leq K\{1 + \mathbf{e}^\lambda(\xi)\}^{1/2} e^{-\lambda s/4},$$

due to the bound (7.12), for some constant $K > 0$. Therefore, by letting $s \to t$ and invoking Lemma 3.3, inequality (7.3) follows. $\square$

**8. Examples.** To illustrate the application of the stated theorems, let us specialize $A = (\Delta - 1)$ in equation (4.1) to get

$$\partial_t^2 u(x,t) = (\Delta - 1)u - 2\alpha\, \partial_t u + f(u,x,t) + \sigma(u, Du, x, t)\, \partial_t W(x,t),$$
$$0 < t < T,\ x \in \mathcal{D} \subset \mathbf{R}^d,\ d \leq 3.$$
(8.1)
$$u(x,0) = u_0(x) \in H^1 \cap L^{2n}, \qquad \partial_t u(x,0) = v_0(x) \in H,$$
$$u(\cdot,t)|_{\partial D} = 0,$$

where $n \geq 1$. We shall present two examples with different kinds of nonlinear terms.

EXAMPLE 1. Let
$$f(u,x,t) = -\kappa u^{2n-1} + \beta(x,t) u^m,$$
(8.2)
$$\sigma(u, Du, x, t) = \zeta(x,t)(1 + |Du|^2)^\delta u^k,$$

where $\kappa > 0$ is a constant, and $\beta$ and $\zeta$ are some functions on $\mathcal{D} \times \mathbf{R}^+$ as yet to be specified. The positive integers $n$, $m$ and $r$ are given such that $2 \leq m < n$, where $n$ is any natural number for $d \leq 2$ but, for $d = 3$, $n \leq 2$. This is so because, according to a Sobolev lemma (see [3], Lemma 4.2), an $L^p$-norm, depending on $p$ and $d$, can be dominated by an $H^1$-norm. Also we assume $\delta \in (0, \frac{1}{2})$ and $0 < k < m(1 - 2\delta)$.

Owing to the Sobolev lemma mentioned above, as in [3], we can show that Conditions A are satisfied so that the equation has a unique global solution. In view of (8.2) and condition (A2), we set $\Phi'(u) = 2\kappa u^{2n-1}$ so that $\Phi(u) = \frac{\kappa}{n}\|u^n\|^2$ and $p_t(u) = \beta(\cdot,t)u^m$. Therefore, we have

(8.3) $\quad (\Phi'(u), u) = 2\kappa \|u^n\|^2 \quad \text{and} \quad \|p_t(u)\|^2 = \|\beta(\cdot,t)u^m\|^2.$

By means of an elementary Young inequality ([9], page 61), it can be shown that, for any $\varepsilon > 0$, we have

$$(\beta u^m)^2 \leq \varepsilon u^{2n} + \frac{\beta^{2q}}{q\varepsilon^{q/q'}},$$

where $q = n/(n-m)$ and $q' = n/m$, so that

(8.4) $\quad \|p_t(u)\|^2 = \int_{\mathcal{D}} \beta^2 u^{2m}\, dx \leq \varepsilon \|u^n\|^2 + C_1 \int_{\mathcal{D}} \beta^{2q}\, dx$

for some constant $C_1 > 0$. Next consider the term

$$\|\Sigma_t(u)\|_R^2 = \int_{\mathcal{D}} r(x,x)\zeta^2(x,t)(1+|Du|^2)^{2\delta} u^{2k}\, dx.$$

By a repeated application of the Young inequality, we can deduce that, for any $\varepsilon', \varepsilon'' > 0$, there exists $C_2 > 0$ such that

(8.5) $\quad \|\Sigma_t(u)\|_R^2 \leq \int_{\mathcal{D}} r(x,x)\zeta^2(x,t)\{\varepsilon' u^{2n} + \varepsilon''|Du|^2 + C_2\}\, dx.$

Suppose that $\beta, \zeta$ and $r$ are bounded and continuous such that

(8.6) $\quad |\beta(x,t)| \leq \beta_0, |\zeta(x,t)| \leq \zeta_0 \quad \text{and} \quad |r(x,x)| \leq r_0$

for any $x \in \mathcal{D}, t \geq 0$. In view of (8.6), inequalities (8.4) and (8.5) yield

(8.7)
$$\|p_t(u)\|^2 \leq \frac{\varepsilon\kappa}{n}\Phi(u) + C_1 \int_{\mathcal{D}} \beta^{2q}(x,t)\, dx,$$
$$\|\Sigma_t(u)\|_R^2 \leq r_0\zeta_0^2\left\{\frac{\varepsilon'\kappa}{n}\Phi(u) + \varepsilon''\|u\|_1^2\right\} + C_2 r_0 \int_{\mathcal{D}} \zeta^2(x,t)\, dx.$$

From (8.3) and (8.4), using the notation in condition (B2), we see that

(8.8)
$$\begin{aligned}\beta_1 &= 2\kappa, & \gamma_1 &= \delta_1 = 0,\\ \beta_2 &= \frac{\varepsilon\kappa}{n}, & \gamma_2 &= 0,\\ \beta_3 &= r_0\zeta_0^2\frac{\varepsilon'\kappa}{n}, & \gamma_3 &= r_0\zeta_0^2\varepsilon''\end{aligned}$$

and

(8.9) $\quad \theta(t) = C_1 \int_{\mathcal{D}} \beta^{2q}(x,t)\, dx, \qquad \rho(t) = C_2 r_0 \int_{\mathcal{D}} \zeta^2(x,t)\, dx.$

Therefore, condition (B3) takes the form

(8.10)
$$2\kappa\lambda^2 - r_0\zeta_0^2\frac{\varepsilon'\kappa}{n}\lambda - \frac{\varepsilon\kappa}{n} > \frac{\lambda^2}{2},$$
$$r_0\zeta_0^2\varepsilon''\lambda < \frac{\lambda^2}{2}.$$

Since $\varepsilon$, $\varepsilon'$ and $\varepsilon''$ are arbitrary, they can be chosen so small that condition (B3) holds simply for $\kappa > \frac{1}{4}$. Assume this is the case. Then, by applying Theorems 4.1, 5.1 and 5.2, depending on the properties of $\theta$ and $\rho$ as defined in (8.9), we can draw the following conclusions:

1. By the conditions in (8.6), it is clear that both $\theta$ and $\rho$ are bounded on $\mathbf{R}^+$ so that, by invoking Theorem 5.1, we can conclude that the solution of the problem (8.1) is bounded in mean-square and there exists a constant $K_1 > 0$ such that
$$\sup_{t>0} E\{\|u_t\|_1^2 + \|\partial_t u_t\|^2 + \|u^n\|^2\} \leq K_1.$$

2. Suppose that the functions $\theta$ and $\rho$ defined by (8.9) belong to $L^1(\mathbf{R}^+)$ so that
$$\int_0^\infty \int_\mathcal{D} \beta^{2q}(x,t)\,dx\,dt < \infty, \qquad \int_0^\infty \int_\mathcal{D} \zeta^2(x,t)\,dx\,dt < \infty.$$

Then, by Theorem 5.2, the solution is ultimately bounded in mean-square if $\kappa > \frac{1}{4}$ and, furthermore, there is $K_2 > 0$ such that
$$E\sup_{t>0}\{\|u_t\|_1^2 + \|\partial_t u_t\|^2 + \|u^n\|^2\} \leq K_2.$$

3. Note that (8.1) has a null solution $(u;v) = (0;0)$. Assume there exists $\alpha_0 > 0$ such that $(e^{\alpha_0 t}\theta)$ and $(e^{\alpha_0 t}\rho)$ belong to $L^1(\mathbf{R}^+)$ or
$$\int_0^\infty \int_\mathcal{D} e^{\alpha_0 t}\beta^{2q}(x,t)\,dx\,dt < \infty, \qquad \int_0^\infty \int_\mathcal{D} e^{\alpha_0 t}\zeta^2(x,t)\,dx\,dt < \infty.$$

Then Theorem 6.1 shows that the null solution is exponentially stable in mean-square and, moreover, according to Theorem 6.3, the solution is in fact a.s. exponentially stable.

EXAMPLE 2. Consider a mildly nonlinear equation of the form

(8.11)
$$\partial_t^2 u(x,t) = (\Delta - 1)u - 2\alpha\,\partial_t u + f(u)$$
$$+ \sigma(Du)\,\partial_t W(x,t), \qquad 0 < t < T,\ x \in \mathcal{D} \subset \mathbf{R}^d,$$
$$u(\cdot,t)|_{\partial D} = 0,$$

subject to the initial conditions $u(x,0) = u_0(x) \in H^1$ and $\partial_t u(x,0) = v_0(x) \in H$, where

$$f(u) = -\kappa u \tan^{-1}(1+u^2),$$
(8.12)
$$\sigma(Du)\partial_t W = \sigma_1\{1+|Du|^2\}^{1/2}\partial_t W^{(1)} + \sigma_2 \partial_t W^{(2)}.$$

In the above equations $\kappa > 0, \sigma_1$ and $\sigma_2$ are some constants, and $W^{(1)}(\cdot,t)$, $W^{(2)}(\cdot,t)$ are independent Wiener random fields with bounded, continuous covariant functions $r_1, r_2$, respectively. Rewriting (8.11) in the system form (7.1) and noting (8.12), it is easy to verify that

(8.13) $$\|F(u)\|^2 \leq \left(\frac{\kappa\pi}{2\eta_1}\right)^2 \|u\|_1^2$$

and

(8.14)
$$\|\Sigma(u)\|_R^2 = \sigma_1^2 \int_{\mathcal{D}} r_1(x,x)(1+|Du|^2)\,dx + \sigma_2^2 \int_{\mathcal{D}} r_2(x,x)\,dx$$
$$\leq \sigma_1^2 r_0 \|u\|_1^2 + c_2$$

for some $c_2 > 0$, where $\eta_1$ is the smallest eigenvalue of $A = (-\Delta + 1)$ and $r_0 = \sup_{x \in \mathcal{D}} |r_1(x,x)|$. Similarly we can obtain the bounds

(8.15)
$$\|F(u) - F(u')\|^2 \leq \frac{\kappa^2}{\eta_1}\left(1+\frac{\pi}{2}\right)^2 \|u-u'\|_1^2,$$
$$\|\Sigma(Du) - \Sigma(Du')\|_R^2 \leq \sigma_1^2 r_0 \|u-u'\|_1^2$$

for any $u, u' \in H^1$. In the notation of Conditions C, we can read off from (8.13) to (8.15) and find $c_1 = 0$,

$$b_1 = \left(\frac{\kappa\pi}{2\eta_1}\right)^2, \qquad b_2 = \sigma_1^2 r_0,$$
$$k_1 = \frac{\kappa^2}{\eta_1}\left(1+\frac{\pi}{2}\right)^2, \qquad k_2 = \sigma_1^2 r_0.$$

To satisfy condition (C3), we require that

$$\left(\frac{\kappa\pi}{2\eta_1}\right)^2 + \sigma_1^2 r_0 \lambda \leq \frac{\lambda^2}{2}$$

and

$$\frac{\kappa^2}{\eta_1}\left(1+\frac{\pi}{2}\right)^2 + \sigma_1^2 r_0 \lambda \leq \frac{\lambda^2}{2}.$$

Then Theorem 7.1 (see the remark following Theorem 5.1) ensures the existence of a unique invariant measure $\mu$ in the state space $\mathcal{H}$ for (8.11) and the the corresponding transition probability converges weakly to $\mu$ at an exponential rate.

## APPENDIX

Let $W(x,t)$ be a continuous Wiener random field as given in Section 3. Then it may be regarded as an $H$-valued Wiener process with a finite-trace covariance operator $R$ with kernel $r(x,y)$. We first define the stochastic integral with an a.s. bounded integrand. To this end, let $\sigma(x,t)$ be an a.s. bounded, continuous predictable random field such that

$$\text{(A.1)} \qquad E \int_0^T \|\sigma_t\|^2 \, dt < \infty.$$

We may consider $\sigma_t$ as a linear operator in $H$ such that $[\sigma_t h](x) = \sigma(x,t)h(x)$ for any $h \in H$. Then it is easy to check that $\sigma_t : H \to H$ is Hilbert–Schmidt a.s. and, noting (A.1),

$$E \int_0^T \text{Tr}(\sigma_t R \sigma_t^\star) \, dt = E \int_0^T \int_{\mathcal{D}} r(x,x)\sigma^2(x,t) \, dx \, dt$$
$$\leq r_0 E \int_0^T \|\sigma_t^2\|^2 \, dt < \infty,$$

where $\star$ denotes the conjugation. Therefore, the stochastic integral

$$\text{(A.2)} \qquad M(x,t) = \int_0^t \sigma(x,s) W(x,ds)$$

or

$$M_t = \int_0^t \sigma_s \, dW_s$$

is well defined as a continuous $H$-valued martingale with mean zero and covariation operator $Q_t$ defined as (see [6], page 90)

$$\text{(A.3)} \qquad \langle (M_\cdot, g), (M_\cdot, h) \rangle_t = \int_0^t (Q_\tau g, h) \, d\tau,$$

where $Q_s$ has the kernel $q(x,y,s) = r(x,y)\sigma(x,s)\sigma(y,s)$.

Now we shall define a stochastic integral with an $L^p$-bounded integrand as shown in the proof of the following theorem.

THEOREM A.1. *Let $W(\cdot, t)$ be a continuous Wiener random field with a bounded covariance function $r(x,y)$ such that*

$$\text{(A.4)} \qquad \sup_{x \in \mathcal{D}} r(x,x) \leq r_0.$$

*Suppose that $\sigma_t = \sigma(\cdot, t)$ is a predictable, continuous $H$-valued process satisfying the condition*

$$\text{(A.5)} \qquad E \int_0^T \|\sigma(\cdot,t)\|^p \, dt = E \int_0^T \int_{\mathcal{D}} |\sigma(x,t)|^p \, dx \, dt < \infty$$

for an integer $p \geq 2$. Then the stochastic integral $M_t$ in (A.2) is well defined as a continuous $H$-valued, $L^p$-martingale with mean zero and covariation operator $Q_t$ for $t \in [0, T]$, as given by (A.3).

PROOF. Since the set $\mathcal{C}_b$ of bounded continuous functions on $\mathcal{D}$ is dense in $L^p(\mathcal{D})$ ([1], page 28), by smoothing, there exists a sequence $\{\sigma_t^n\}$ of predictable continuous random fields converging to $\sigma_t$ such that it satisfies condition (A.1) and

$$\text{(A.6)} \qquad \lim_{n \to \infty} E \int_0^T \|\sigma_t^n - \sigma_t\|^p \, dt = 0.$$

Therefore, as in (A.2), the stochastic integral

$$M_t^n = \int_0^t \sigma_s^n \, dW_s$$

exists as a continuous $H$-valued martingale for each $n$. Let $\mathcal{M}_T^p$ denote the Banach space of continuous $L^p$-martingales $N_t \in H$ with norm ([6], page 79)

$$\|N\|_T = \left\{ E \sup_{0 \leq t \leq T} \|N_t\|^p \right\}^{1/p}.$$

Then the sequence $\{M_t^n\}$ belongs to $\mathcal{M}_T^p$, since, by the B–D–G inequality,

$$\|M^n\|_T^p = E \sup_{0 \leq t \leq T} \|M_t^n\|^p$$

$$\leq C_p E \left\{ \int_0^T \int_{\mathcal{D}} r(x,x) |\sigma^n(x,t)|^2 \, dx \, dt \right\}^{p/2}$$

$$\leq C_p r_0^{p/2} E \left\{ \int_0^T \|\sigma_t^n\|^2 \, dt \right\}^{p/2}$$

$$\leq C_p(T) E \int_0^T \|\sigma_t^n\|^p \, dt,$$

in which, in view of (A.5), the upper limit is bounded, where $C_p, C_p(T)$ are some positive constants.

Now, for $n > m$,

$$\|M^n - M^m\|_T^p = E \sup_{0 \leq t \leq T} \|M_t^n - M_t^m\|^p$$

$$\leq C_p E \left\{ \int_0^T \text{Tr} \, Q_s^{mn} \, ds \right\}^{p/2},$$

where

$$\text{Tr} \, Q_s^{mn} = \int_D q^{mn}(x,x,s) \, dx$$

$$= \int_D r(x,x)[\sigma^n(x,s) - \sigma^m(x,s)]^2 \, dx$$
$$\leq r_0 \|\sigma_s^n - \sigma_s^m\|^2.$$

It follows from (A.4) and (A.6) that

$$\|M^n - M^m\|_T^p \leq C_p(T) E \int_0^T \|\sigma_t^n - \sigma_t^m\|^p \, dt,$$

which goes to zero as $n > m \to \infty$ due to (A.6). Therefore, the sequence $\{M_t^n\}$ converges to the limit denoted by $M_t$, which is defined as a stochastic integral given by (A.2). We can check that it preserves the properties of $M_t^n$ as stated in the theorem. $\square$

**Acknowledgments.** The author wishes to thank the referees for their thorough reviews of the original manuscript and helpful suggestions which led to a substantial improvement in the presentation. He is indebted to Professor Jerzy Zabczyk of the Polish Academy of Science for a helpful discussion about the problem of invariant measures given in Theorem 7.1.

Department of Mathematics
Wayne State University
Detroit, Michigan 48202
USA
E-mail: plchow@math.wayne.edu